\documentclass[10pt]{amsart}

\usepackage{amsmath,amssymb,amscd,amsxtra}
\usepackage{graphics,graphicx}
\usepackage{latexsym}

\newcommand{\N}{\mathbb{N}}
\newcommand{\R}{\mathbb{R}}
\newcommand{\C}{\mathbb{C}}
\newcommand{\Z}{\mathbb{Z}}
\newcommand{\Q}{\mathbb{Q}}

\renewcommand{\S}{\mathcal{S}}

\newcommand{\ip}[2]{\langle#1,#2\rangle}

\newcommand{\ga}[3]{M_{{#1}\beta }T_{{#2}\alpha}{#3}}

\newcommand{\cG}{\mathcal{G}}

\newtheorem{prop}{Proposition}

\theoremstyle{remark}

\theoremstyle{definition}

\theoremstyle{conjecture}
\newtheorem{conjecture}{Conjecture}
\theoremstyle{claim}

\newtheorem{question}{Question}

\begin{document}

\title{An invitation to Gabor analysis}

\author{Kasso A.~Okoudjou}

\address{Kasso A.~Okoudjou\\
Department of Mathematics, Massachusetts Institute of Technology\\
Cambridge, MA 02139  USA}

\email{kasso@mit.edu}
\thanks{This work  was partially supported by a grant from the Simons Foundation $\# 319197$, ARO grant W911NF1610008,  the National Science Foundation under Grant No.~DMS-1814253, and an MLK  visiting professorship}

\subjclass[2000]{Primary 42C15, 42C4; Secondary 46C05}

\date{\today}

\keywords{HRT conjecture, Gabor Analysis,  short-time Fourier transform, time-frequency analysis, Wilson bases}

\begin{abstract} Gabor analysis which can be traced back to Dennis Gabor's influential 1946 paper \emph{Theory of communication},  is concerned with both the theory and the applications of  the approximation properties  of sets of time and frequency shifts of a given window function. It re-emerged with the advent of wavelets at the end of the last century and is now at the intersection of many fields of mathematics, applied mathematics, engineering, and science.  The goal of this paper is to give a brief introduction to Gabor analysis by elaborating on three open problems.
\end{abstract}

\maketitle \pagestyle{myheadings} \thispagestyle{plain}
\markboth{K. A. OKOUDJOU}{AN INVITATION TO GABOR ANALYSIS}

\section{Introduction}\label{sec1}
Using the ubiquitous theory of Fourier series,  one can decompose and reconstruct  any $1$-periodic and square integrable function in terms of complex exponential functions with frequencies at the integers. More specifically, for any such function $f$ we have $$f(x)=\sum_{n=-\infty}^{\infty}c_n e^{2\pi i nt}$$ where the coefficients $\{c_n\}_{n\in \Z}$ are square summable and the series converges in mean square, that is $$\lim_{N\to \infty}\int_0^1\bigg| f(t)-\sum_{n=-N}^Nc_ne^{2\pi i nt}\bigg|^2 dt =0.$$ The significance of this simple fact is that $f$ is completely determined by the coefficients $\{c_n\}$, and, conversely, each square summable sequence gives rise to a  unique $1$-periodic and square integrable function. This fact is equivalent to saying that  the set  $\{e_n(t):=e^{2\pi i nt}\}_{n=-\infty}^{\infty}$  forms an orthonormal basis (ONB) for $L^2([0,1))$. We shall consider these functions  as the building blocks of Fourier analysis of the space of $1$-periodic square integrable functions.

 In his celebrated work \cite{gabo1946}, Dennis Gabor sought to decompose any square integrable function on the real line in a similar manner. To this end, he proposed to ``localize'' the Fourier series decomposition of such a function, by  first using translates of  an appropriate window function to restrict the function to  time intervals that cover the real line. The next step in the process is to write the Fourier series of each of the ``localized functions'', and finally, one superimposes all these local Fourier series.  Putting this into practice, Gabor chose the Gaussian as a window and claimed that every square integrable function $f$ on $\R$   has the following (non orthogonal) expansion 

\begin{equation}\label{gauexp}
f(x)=\sum_{n\in \Z}\sum_{k\in \Z}c_{nk}\, e^{-\tfrac{\pi (x-n\alpha)^2}{2\alpha^2}}\, e^{2\pi i kx/\alpha}
\end{equation} where $\alpha>0$. Furthermore, he  argued on how to find the coefficients $(c_{nk})_{n, k \in \Z}\in \C$ using successive local approximations by Fourier series.  In fact,  in 1932,  John von Neumann already made a related claim, when he stipulated  that the system of functions 
\begin{equation}\label{gabgauss}
\mathcal{G}(\varphi, 1,1)=\big\{\varphi_{ nk}(\cdot):=e^{2\pi i k \cdot} \varphi(\cdot -n): n, k\in \Z\big\}
\end{equation} where $\varphi(x)=e^{-\pi x^2}$ spans a dense subspace of $L^2(\R)$ \cite{JvNeu}.

Both claims were  positively established in $1971$ independently by V.~Bargmann, P.~Butera, L.~Girardelo, and J.~R.~Klauder \cite{Bargmann71} and A.~M.~Perelomov \cite{Perelomov71}. 
 It then follows that both statements hint at the fact that any square integrable function $f$ is completely determined in the time-frequency plane by  the coefficients $\{c_{n k}\}_{k, n \in \Z}$.    In contrast to the theory of Fourier series, the building blocks in this process are the time and frequency shifts of a function such as the Gaussian: $\{\varphi_{nk}(x)=e^{2\pi i k \cdot} \varphi(\cdot -n): n, k\in \Z\}$. But as we shall see later, we could consider time-frequency shifts of other square integrable functions along a lattice $\alpha \Z \times \beta\Z$ leading to $\{e^{2\pi i \beta k \cdot}g(\cdot - \alpha n): k, n\in \Z\}$. The main point here is that the building blocks can dependent on three parameters: $\alpha>0$  corresponding to shifts in time/space, $\beta>0$ representing shifts in frequency, and a square integrable window function $g$. 

In some sense, both Gabor and von Neumann statements  can also be  thought of as the foundations of what is known today as Gabor analysis, an active research field at the intersection of (quantum) physics, signal processing, mathematics, and  engineering. In broad terms, Gabor analysis seeks to develop  (discrete) joint  time/space-frequency representations  of  functions (distributions, or signals) initially defined only in time or frequency, and it re-emerged with the advent of wavelets \cite{DaGroMe, Dau90, Daub92, HeiWal}.  For a more complete introduction to the theory and applications of Gabor analysis we refer to  \cite{chri2003, FeiStr1, FeiStr2, Groc2001, Heil07, WeissHern}.

The goal of this paper is to give an overview of some  interesting  open problems in Gabor analysis that are in need of solutions. But first, in Section~\ref{sec2} we review  some fundamental results in Gabor analysis. In Section~\ref{sec3}, we consider the problem of characterizing the set of all ``good'' parameters $\alpha, \beta$ for a fixed window function $g$. In Section~\ref{sec4} we consider the problem of constructing orthonormal bases for $L^2(\R)$ by taking appropriate (finitely many) linear combinations of time-frequency shifts of $g$ along a lattice $\alpha \Z \times \beta\Z$. Finally, in Section~\ref{sec5} we elaborate on a conjecture that asks whether any finite set of time-frequency shifts of a square integrable function is linearly independent.

\section{Gabor frame theory}\label{sec2} 
We start with a motivating example based on the  $L^2$  theory of Fourier series. In particular, we would like to exhibit a set of building blocks $\{g_{nk}\}_{k, n\in \Z}$ that can be used to decompose every square integrable function. To this end, let $g(x)=\chi_{[0,1)}(x)$, where $\chi_{I}$ denotes the indicator function of the measurable set  $I.$ 
Any $f\in L^2(\R)$ can be localized to the interval $[n, n+1)$ by considering its restriction, $f(\cdot)\, g(\cdot -n) $ to this interval. By superimposing all these restrictions over all integers $n\in \Z$, we recover the function $f$. That is, we can write 
\begin{equation}\label{eq:partition}
f(x)=\sum_{n=-\infty}^\infty f(x)g(x-n)
\end{equation}
 with convergence $L^2$.  But since the restriction of $f$ to $[n, n+1)$ is square integrable,  it  can be expanded into its $L^2$ convergent Fourier series leading to 
\begin{equation}\label{eq:locfs}
f(x)\, g(x -n)=\sum_{k\in \Z}c_{nk} e^{2\pi i x k}
\end{equation}
 where for each $k\in \Z,$ 
  $$c_{nk}=\ip{f(\cdot) g(\cdot -n)}{e^{2\pi i k\cdot}}_{L^2([n, n+1))}=\int_{-\infty}^{\infty}f(x)\, g(x -n)\, e^{-2\pi i kx}\, dx=\ip{f}{g_{nk}}$$ 
with $g_{nk}(x)=g(x -n)\, e^{2\pi i kx}$. Here and in the sequel, $\ip{\cdot}{\cdot}$ denotes the inner product on either $L^2(\R)$, the space Lebesgue measurable square integrable functions on $\R$, or $\ell^2(\Z^2)$ the space of square summable sequences on $\Z^2$.  In addition, we use the notation $\| \cdot\|:= \|\cdot\|_2$ to denote the corresponding norm. The context will make it clear which of the two spaces we are dealing with.  

Substituting this in~\eqref{eq:locfs} and~\eqref{eq:partition} leads to 
\begin{equation}\label{eq:chigab}
f(x)=\sum_{n=-\infty}^\infty f(x)g(x-n)=\sum_{n=-\infty}^{\infty}\sum_{k=-\infty}^{\infty}c_{nk} e^{2\pi i kx}g(x-n)= \sum_{k, n=-\infty}^{\infty}\ip{f}{g_{nk}}g_{nk}(x)
\end{equation}

This expansion of $f$  is similar to Gabor's claim~\eqref{gauexp}, with the the following key differences:\newline
\noindent $\bullet$ The coefficients in~\eqref{eq:chigab} are explicitly given and are linear in $f$. \newline
\noindent $\bullet$ ~\eqref{gauexp} is based on the Gaussian while~\eqref{eq:chigab} is based on the indicator function of $[0,1)$. \newline
\noindent $\bullet$ Finally, the expansion given in~\eqref{eq:chigab} is is an orthonormal decomposition while the one given by~\eqref{gauexp} is not. 

One of the goals of this section is to elucidate the difference in behavior between the two building blocks appearing in~\eqref{gauexp} and~\eqref{eq:chigab}. In addition, we shall elaborate on the existence of orthonormal bases of the form $\{e^{2\pi i \beta k \cdot}g(\cdot - \alpha n): k, n\in \Z\}$.

The two systems of functions in~\eqref{gauexp} and~\eqref{gabgauss} are examples of Gabor (or Weyl-Heisenberg) systems. More specifically, for $a, b \in \R$  and a function $g$ defined on $\R$, let $M_bf(x)=e^{2\pi i bx}f(x)$ and $T_af(x)=f(x-a)$  be  respectively  the modulation operator, and the translation operator.   The \emph{Gabor system} generated by a function $g\in L^{2}(\R)$, and parameters $\alpha, \beta>0$,  is the set of functions \cite{Groc2001}  $$\mathcal{G}(g, \alpha, \beta )=\{M_{k\beta}T_{n\alpha}g(\cdot)=e^{2\pi i k \beta  \cdot}g(\cdot - n\alpha  ): k, n \in \Z\}.$$

Given $g\in L^{2}(\R)$, and $\alpha, \beta>0$,  the Gabor system  $\mathcal{G}(g, \alpha, \beta )$ is called a \emph{frame} for $L^2(\R)$ is there exist constants $0<A\leq B< \infty$ such that 
\begin{equation}\label{gabframe}
A\|f\|^2\leq \sum_{k, n\in \Z}|\ip{f}{M_{k\beta}T_{n\alpha}g}|^2\leq B\|f\|^2 \qquad \forall f\in L^2(\R).
\end{equation} The constant $A$ is called a lower frame bound, while  $B$ is called an upper frame bound.  When $A=B$ we  say that the Gabor frame is \emph{tight}.  In this case, the frame bound $A$ is referred to as the \emph{redundancy} of the frame. Loosely speaking, the redundancy $A$ measures by how much the Gabor tight frame is overcomplete. 
A tight Gabor frame for which $A=B=1$ is called a \emph{Parseval} frame. Clearly, if $\mathcal{G}(g, \alpha, \beta )$ is an ONB then it is a Parseval frame, and conversely,  if $\mathcal{G}(g, \alpha, \beta)$ is a Parseval frame and $\|g\|=1$, then it is a Gabor ONB.  

More generally, a Gabor frame is a ``basis-like'' system that can be used to decompose and/or reconstruct any square integrable function. As such, it will not come as a surprise that generalizations of certain tools from linear algebra might be useful in analyzing Gabor frames. 
We refer to \cite{chri2003, Daub92, Groc2001, WeissHern, Heil07}  for more background on Gabor frames, and  summarize below some results needed in the sequel. 

Suppose we would like to analyze $f$ using the Gabor system $\mathcal{G}(g, \alpha, \beta).$ We are then led to consider the correspondence that takes any square integrable function $f$ into the sequence $\{\ip{f}{ M_{k\beta }T_{n\alpha }g}\}_{k, n\in \Z}$. This correspondence is sometimes called the \emph{analysis or decomposition operator} and denoted by $$C_g: f\to \{\ip{f}{ M_{k\beta }T_{n\alpha }g}\}_{k, n\in \Z}.$$ Its (formal) adjoint $C_g^{*}$ called the \emph{synthesis or reconstruction operator} maps sequences $c=\{c_{kn}\}_{k, n\in \Z}$ to $$C^{*}_gc=\sum_{k, n\in \Z}c_{kn}M_{k\beta }T_{n\alpha }g.$$ The composition of these two operators is called the \emph{(Gabor) frame operator} associated to the Gabor system $\mathcal{G}(g, \alpha, \beta)$ is defined by 
\begin{equation}\label{gafop}
Sf:=S_{g, \alpha, \beta}f=C_g^{*}C_g (f)=\sum_{n, k\in \Z}\ip{f}{M_{k\beta }T_{n\alpha}g}M_{k\beta }T_{n\alpha}g
\end{equation}

It follows that, given $f\in L^2(\R)$, we can (formally) write that $$\ip{Sf}{f}=\ip{C_{g}^{*}C_gf}{f}=\ip{C_gf}{C_gf}=\sum_{k, n\in \Z}|\ip{f}{M_{k\beta }T_{n\alpha }g}|^2.$$   Therefore,   $\cG(g, \alpha, \beta)$ is a frame for $L^2$ if and only if there exist constants $0<A\leq B< \infty$ such that $$A\|f\|_2^2\leq \ip{Sf}{f}\leq B\|f\|^2_2\qquad \forall\, f\in L^2(\R).$$ In particular, $\cG(g, \alpha, \beta)$ is a frame for $L^2$ if and only if the self-adjoint frame operator $S$ is bounded and positive definite. Furthermore, the optimal upper frame bound $B$ is the largest eigenvalue of $S$ while the optimal lower bound $A$ is its smallest eigenvalue.  In addition, $\cG(g, \alpha, \beta)$ is a tight  frame for $L^2$ if and only if $S$ is a multiple of the identity. 
   
Viewing a Gabor frame as an overcomplete  ``basis-like'' object suggests that  any square integrable function can be written in a non-unique way as a linear combination of the Gabor atoms $\{M_{k\beta }T_{n\alpha }g\}_{k, n \in \Z}$.  Akin to  the role of the pseudo-inverse in linear algebra, we single out one expansion that results in a somehow canonical representation of $f$ as a linear combination of $\{M_{k\beta }T_{n\alpha }g\}_{k, n \in \Z}$.  To obtain this decomposition we need a few basic facts about the frame operator.

Suppose that $\cG(g, \alpha, \beta)$ is a Gabor frame for $L^2$, and let $f\in L^2$. For all $(\ell, m) \in \Z^2$ the frame operator $S$ and $M_{\ell \beta}T_{m\alpha}$ commute. That is  $$S(M_{\ell \beta}T_{m\alpha}f)=M_{\ell \beta}T_{m\alpha}(S(f)) \, \textrm{for\, all\,} \, (\ell, m) \in \Z^2.$$ It follows that $S^{-1}$  and $M_{\ell \beta}T_{m\alpha}$ also commute  for all $(\ell, m) \in \Z^2$. As a  consequence, given  $f\in L^2(\R)$ we have
\begin{align*}
f&=S(S^{-1}f)=\sum_{k, n\in \Z}\ip{S^{-1}f}{\ga{k}{n}{g}}\ga{k}{n}{g}\\
&= \sum_{k, n\in \Z}\ip{f}{S^{-1}\ga{k}{n}{g}}\ga{k}{n}{g}=\sum_{k, n\in \Z}\ip{f}{\ga{k}{n}{\tilde{g}}}\ga{k}{n}{g}\end{align*}
 where $\tilde{g}=S^{-1}g\in L^2(\R)$ is called the \emph{canonical dual} of $g$. Similarly, by writing $f=S^{-1}(Sf)$ we get that $$f=\sum_{k, n\in \Z}\ip{f}{\ga{k}{n}{g}}\ga{k}{n}{\tilde{g}}.$$
 
 The coefficients $\{\ip{f}{\ga{k}{n}{\tilde{g}}}\}_{k, n\in \Z}$ give the least square approximation of $f$. Indeed,  for  $f \in L^2$, let $\tilde{c}=(\ip{f}{\ga{k}{n}{\tilde{g}}})_{k, n\in \Z}\in \ell^2(\Z^2).$ Given any (other) sequence $(c_{k, n})_{k, n\in \Z} \in \ell^2(\Z^2)$  such that $$f=\sum_{k, n\in \Z}\tilde{c}_{k,n}\ga{n}{k}{g}=\sum_{k, n\in \Z}c_{k,n}\ga{k}{n}{g},$$  we have 

$$\|\tilde{c}\|_2^2=\sum_{k, n\in \Z}| \ip{f}{\ga{k}{n}{\tilde{g}}}|^2=\ip{S^{-1}f}{f}=\sum_{k,n\in \Z}c_{k,n}\ip{S^{-1}\ga{k}{n}{g}}{f}=\sum_{k,n\in \Z}c_{k,n}\overline{\tilde{c}_{{k,n}}}=\ip{c}{\tilde{c}}.$$
Consequently, $\ip{c-\tilde{c}}{\tilde{c}}=0$, leading to 
 $$\|c\|^2_2=\|c-\tilde{c}\|^2_2+\|\tilde{c}\|^2_2\geq \|\tilde{c}\|^2_2$$ with equality if and only if $c=\tilde{c}.$ In other words, for a Gabor frame $\cG(g, \alpha, \beta)$, and given $f\in L^2$, among all expansions $f= \sum_{k, n\in \Z}c_{k,n}\ga{k}{n}{g},$ with $c=(c_{k,n})_{k, n \in \Z}\in \ell^2(\Z^2)$, the coefficient $\tilde{c}= (\ip{f}{\ga{k}{n}{\tilde{g}}})_{k, n\in \Z}\in \ell^2(\Z^2)$ has the least norm.

  Because the frame operator $S$ is positive definite, $S^{1/2}$ is well defined and positive definite as well. Thus, we can write $$f=S^{-1/2}SS^{-1/2}f= \sum_{k, n}\ip{f}{S^{-1/2}\ga{k}{n}{g}}S^{-1/2}\ga{k}{n}{g}= \sum_{k, n}\ip{f}{\ga{k}{n}{g^{\dag}}}\ga{k}{n}{g^{\dag}}$$ where $g^{\dag}=S^{-1/2}g\in L^2$. In other words, $\mathcal{G}(g^{\dag}, \alpha, \beta)$ is a Parseval frame.

 Finally, assume that $A, B$ are the optimal frame bounds for $\cG(g, \alpha, \beta)$.  Then, for all $f\in L^2$, we have
 
  $$\sum_{k, n\in Z}|\ip{f}{\ga{k}{n}{\tilde{g}}}|^2=\ip{S^{-1}f}{f}=\ip{S^{-1}f}{S(S^{-1}f)}\leq B\|S^{-1}f\|^2\leq \tilde{B}\|f\|^2$$ and similarly, we have the lower bound $$\ip{S^{-1}f}{f}=\ip{S^{-1}f}{S(S^{-1}f)}\geq A\|S^{-1}f\|^2\geq \tilde{A}\|f\|^2$$ Therefore, if $\cG(g, \alpha, \beta)$ is Gabor frame for $L^2(\R)$, then so is $\cG(\tilde{g}, \alpha, \beta)$ where $\tilde{g}=S^{-1}g\in L^2(\R)$.   We summarize all these facts in the following result.

\begin{prop}[Reconstruction formulas for Gabor frame]\label{gaborf_properties} Let $g\in L^2(\R)$ and $\alpha, \beta>0$. Suppose that $\cG(g, \alpha, \beta)$ is a frame for $L^2(\R)$ with frame bounds $A, B$. Then the following statements hold.
\begin{enumerate}
\item[(a)] The Gabor system $\mathcal{G}(\tilde{g}, \alpha, \beta)$ with $\tilde{g}=S^{-1}g\in L^2$, is also a frame for $L^2$ with frame bounds $1/B, 1/A$. Furthermore, for each $f\in L^2$ we have the following reconstruction formulas:

$$f=\sum_{k, n\in \Z}\ip{f}{\ga{k}{n}{\tilde{g}}}\ga{k}{n}{g}=\sum_{k, n\in \Z}\ip{f}{\ga{k}{n}{g}}\ga{k}{n}{\tilde{g}}.$$  In addition, among all sequences $c=(c_{k, n})_{k, n\in \Z} \in \ell^2(\Z^2)$  such that $f=\sum_{k, n\in \Z}c_{k,n}\ga{k}{n}{g}$, the sequence $\tilde{c}=(\ip{f}{\ga{k}{n}{\tilde{g}}})_{k, n\in \Z}\in \ell^2(\Z^2)$ satisfies $$\|\tilde{c}\|^2_2=\sum_{k, n\in \Z}| \ip{f}{\ga{k}{n}{\tilde{g}}}|^2\geq \sum_{k, n\in \Z}|c_{k,n}|^2=\|c\|^2_2$$ with equality if and only if $c=\tilde{c}.$
\item[(b)] The Gabor system $\mathcal{G}(g^{\dag}, \alpha, \beta)$ where $g^{\dag}=S^{-1/2}g\in L^2$,  is a Parseval frame. In particular, each $f\in L^2$ has the following expansion 
$$f= \sum_{k, n}\ip{f}{\ga{k}{n}{g^{\dag}}}\ga{k}{n}{g^{\dag}}.$$ 
\end{enumerate} 
\end{prop}

 It is worth pointing out that the coefficients $(\ip{f}{\ga{k}{n}{g}})_{k, n}$ appearing in~\eqref{gabframe} or in~\eqref{gafop} are samples of the \emph{Short-Time Fourier Transform} (STFT) of $f$ with respect to $g$. This is the function  $V_g$ defined on $ \R^2$ by $$V_gf(x, \xi)=\ip{f}{M_{\xi}T_xg}=\int_{\R}f(t)\overline{g(t-x)}e^{-2\pi i t\xi}\, dt.$$  When $g\in L^2(\R)$ is chosen such that $\|g\|=1$, then $V_g$ is an isometry from $L^2(\R)$ onto a closed subspace of $L^2(\R^2)$ and for all $f\in L^2(\R)$ 
 \begin{equation}\label{stft-iso}
 \int_{\R}|f(t)|^2\, dt = \iint_{\R^2}|V_gf(x, \xi)|^2 \, dx\, d\xi
 \end{equation}
  Furthermore,  for any $h\in L^2$ such that $\ip{g}{h}\neq 0$ 
  \begin{equation}\label{stft-inv}
  f(t)=\tfrac{1}{\ip{g}{h}}\iint_{\R^2}V_gf(x, \xi)M_{\xi}T_xh(t)\, dx\, d\xi
  \end{equation}
   where the integral is interpreted in the weak sense. We refer to \cite[Chapter 1]{FollHAPS89} and  \cite[Chapter 3]{Groc2001} for more on the STFT and related phase-space or time-frequency transformations.

The reconstruction formulas in  Proposition~\ref{gaborf_properties} can be viewed as  discretizations of  the inversion formula for the STFT~\eqref{stft-inv}. In particular, sampling the STFT on the lattice $\alpha \Z \times \beta \Z$ and using the weights $\tilde{c}=(\ip{f}{\ga{k}{n}{\tilde{g}}})_{k, n\in \Z}=(V_{\tilde{g}}f(\alpha k, \beta n))_{k, n \in \Z} \in \ell^2(\Z^2)  $ perfectly reconstructs $f$. As such one can expect that in addition to the quality of the window $g$ (and hence $\tilde{g}$), the density of the lattice must play a role in establishing these formulas. Thus, it must not come as a surprise that the following results hold.

\begin{prop}[Density theorems for Gabor frames]\label{density_gabf} Let $g\in L^2(\R)$ and $\alpha, \beta>0$.

\begin{enumerate}
\item[(a)] If $\cG(g, \alpha, \beta)$ is a Gabor frame for $L^2(\R)$ then $0<\alpha\beta \leq 1$
\item[(b)] If $\alpha\beta>1$, then $\cG(g, \alpha, \beta)$ is incomplete in $L^2(\R)$. 
\item[(c)] $\cG(g, \alpha, \beta)$ is an orthonormal basis for $L^2(\R)$  if and only $\cG(g, \alpha, \beta)$ is a tight frame for $L^2(\R)$, $\|g\|=1$, and $\alpha\beta=1$.
\end{enumerate}
\end{prop}
These results were proved using various techniques ranging from operator theory to signal analysis illustrating the multi-origin of  Gabor  frame theory. For more on these density results we refer to \cite{Bagg90, Dau90, DauLL95, Foll06, RT95}, and for a  complete historical perspective see \cite{Heil07}.

 At this point some questions  arise naturally. For example,  can one classify  $g\in L^2(\R)$ and the parameters $\alpha, \beta>0$, such that  $\mathcal{G}(g, \alpha, \beta )$  generates a frame, or an ONB for $L^2(\R)$? 
  Despite some spectacular results both in the theory and the applications of Gabor frames \cite{FeiStr1, FeiStr2}, these problems have not been completely  resolved.  Section~\ref{sec3} will be devoted to addressing the \emph{frame set problem for Gabor frames}. That is given $g\in L^2(\R)$ characterize the set of all $(\alpha, \beta)\in \R^2_+$ such that $\mathcal{G}(g, \alpha, \beta)$ is a frame.  On the other hand, and as seen from part (c) of Proposition~\ref{density_gabf}, Gabor ONB can only occur when $\alpha \beta=1$. In addition to this restriction,  there does not exist a Gabor ONB with  $g\in L^2$ such that $$\int_{-\infty}^{\infty} |x|^2\, |g(x)|^2\, dx \int_{-\infty}^{\infty}|\xi|^2\, |\hat{g}(\xi)|^2\, d\xi < \infty$$  where $$\hat{g}(\xi)=\int_{-\infty}^{\infty}g(t)e^{-2\pi i t \xi}dt$$ is the Fourier transform of $g$.  This uncertainty principle-type result known as the Balian-Low Theorem (BLT) precludes the existence  of Gabor ONBs with well-localized windows \cite{Bali81, Bat88, Dau90, Low85}. We use the term well-localized window to describe functions $g$ that behave well in both time/space and frequency. For example, functions in certain Sobolev spaces, and more generally in the so-called modulation spaces can be thought of as well-localized \cite[Chapter 11]{Groc2001}. With this in mind, the following result holds. 
  
  \begin{prop}[The Balian-Low Theorem]\label{blt}
  Let $g \in L^2(\R)$ and $\alpha>0$. If $\mathcal{G}(g, \alpha, 1/\alpha )$ is an orthonormal basis for $L^2(\R)$ then 
  $$\int_{-\infty}^{\infty} |x|^2\, |g(x)|^2\, dx\,  \int_{-\infty}^{\infty}|\xi|^2\, |\hat{g}(\xi)|^2\, d\xi =\infty.$$
   \end{prop}
  We refer to \cite{BenHeiWal94} for a survey on the BLT.   In Section~\ref{sec4} we will introduce a modification of Gabor frames that will result in an ONB called \emph{Wilson} basis with well-localized (or regular) window functions $g$. These ONBs were introduced by K.~G.~Wilson \cite{wi87} under the name of Generalized Warnnier functions. The fact that these are indeed ONBs was later established by Daubechies, Jaffard and Journ\'e \cite{DaJaJo91} who developed  a systematic construction method for these kinds of systems. The method starts with constructing a tight Gabor frame of redundancy $A=2$ and a well-localized window $g$. By then taking appropriate linear combinations of at most two Gabor atoms from this tight Gabor frame, the author removed the original redundancy and obtained an ONB. While it is clear that tight Gabor frame with well-localized generators and arbitrary redundancy can be constructed, it remains an open question how or if one can get ONBs from these systems. We survey this question in Section~\ref{sec4}, and mention that an interesting application involving the Wilson bases is the recent detection of the gravitational waves  \cite{NeKlMi12, Drag16}. We refer to the short survey \cite{ChMJM16}, and to \cite{Mey18} for some historical perspectives. For more on the Wilson bases we refer to \cite{BJLO17, KutStr,  Woj07, Woj11}.

 \section{The frame set problem for Gabor frames}\label{sec3}
 As mentioned in the Introduction, a Gabor system is determined by three parameters: the shift parameters $\alpha, \beta$, and the window function $g$. Ideally, one would like to classify the set of all these three parameters for which the resulting system is a frame. However, and in general, this is a difficult question and we shall only consider the special case in which the window function $g$ is fixed and one seeks the set of all parameters $\alpha, \beta>0$ for which the resulting system is a frame.

 In this setting, the \emph{frame set} of a  function $g\in L^2(\R)$  is defined as $$\mathcal{F}(g)=\left\{(\alpha, \beta)\in \R_+^2:\, \cG(g,\alpha, \beta)\ \mbox{is a frame}\right\}
.$$  In general, determining $\mathcal{F}(g)$ for a given function $g$ is also an open problem. One of the most known general result proved by Feichtinger and Kaiblinger \cite{FeiKai04} states that $\mathcal{F}(g)$ is an open subset of $\R_+^2 $ if $g \in L^2(\R)$ belongs to the modulation space $M^1(\R)$ (\cite{Groc2001}), i.e., $$\iint_{\R^2}|V_gg(x, \xi)|\, dxd\xi<\infty. $$ Examples of functions in this space include $g(x)=e^{-\pi |x|^2}$ or $g(x)=\tfrac{1}{\cosh x}$. In fact, for these specific functions more is know. Indeed, $$\mathcal{F}(g)= \left\{(\alpha, \beta)\in \R_+^2:\,  \alpha \beta<1 \right\}$$ if   $g\in \{e^{-\pi x^2}, \tfrac{1}{\cosh x}, e^{-x}\chi_{[0,+\infty]}(x), e^{-|x|}\}$,  \cite{Jan, Janss, JanStroh, Lyuba, Seip, SeiWal}. 
On the other hand when $g(x)=\chi_{[0,c]}(x), c>0$,  $\mathcal{F}(g)$ is a rather complicated set that has only been fully described in recent years by Dai and Sun \cite{DaiSun}, see also  \cite{GuHan, Jans} for earlier work on this example.

Let $g(x)=e^{-|x|}$ and observe that $\hat{g}(\xi)=\tfrac{2}{1+4\pi^2\xi^2}$, which makes  $g(x)=e^{-|x|}$   an example of a totally positive function of type $2$. More generally,  $g\in L^2(\R)$ is \emph{a totally positive function of type} $M$, where $M$ is a natural number, if its Fourier transform has the form $\hat{g}(\xi)=\prod_{k=1}^M (1+2\pi i \delta_k \xi)^{-1}$ where $\delta_k\neq \delta_\ell \in \R$ for $k\neq \ell$. It was proved that for all such functions $g$,  $\mathcal{F}(g)= \left\{(\alpha,\beta)\in \R_+^2:\,  \alpha \beta<1 \right\}$  \cite{Grosto, Groch}. A similar result holds for the class of totally positive functions of Gaussian type, \cite{GrRoSt}, which are functions whose Fourier transforms have the form $\hat{g}(\xi)=\prod_{k=1}^M (1+2\pi i \delta_k \xi)^{-1} e^{-c\xi^{2}}$ where $\delta_1\neq \delta_2\neq  \hdots \neq \delta_M \in \R$ and $c>0$. We refer to  \cite{Groch} for a survey of the structure of $\mathcal{F}(g)$ not only for the rectangular lattices we consider here, but more general Gabor frame on discrete (countable) sets $\Lambda \subset \R^2$.

However, there are other ``simple'' functions $g$ for which determining $\mathcal{F}(g)$ remains largely a mystery. In the rest of this section we consider the frame set for the $B$ splines  $g_N$ given by 
$$\begin{cases} g_1(x)= \chi_{[-1/2, 1/2]},\quad \text{and}\\
g_N(x)=g_1\ast g_{N-1}(x)\quad \text{for}\quad N\geq 2. \end{cases}$$ 
Christensen lists  the characterization of $\mathcal{F}(g_N)$ for $N\geq 2$ as one  of the six main problems in frame theory \cite{Ole1}. Due to the fact that  $g_N \in M^1(\R)$ for $N\geq 2$, we know that $\mathcal{F}(g_N)$ is an open subset of $\R^2_+$. The current description of points in this set can be found in \cite{AtiKouOko1, Olehon1,  Olehon, KarGro, Kloosto, Lemniel}.

For example, consider the case $N=2$ where $$g_2(x)= \chi_{[-1/2, 1/2]}\ast \chi_{[-1/2, 1/2]}(x)= \max{(1-|x|,0)}=\left\{
 \begin{array}{r@{\quad \text{if} \quad}l}
  1+x &  x\in[-1, 0] \\
  1-x & x\in[0,1]
 \end{array}\right.$$
The known results on $\mathcal{F}(g_2)$ can be summarized as follows.  

\begin{figure}[!h]
 \includegraphics[scale=0.3]{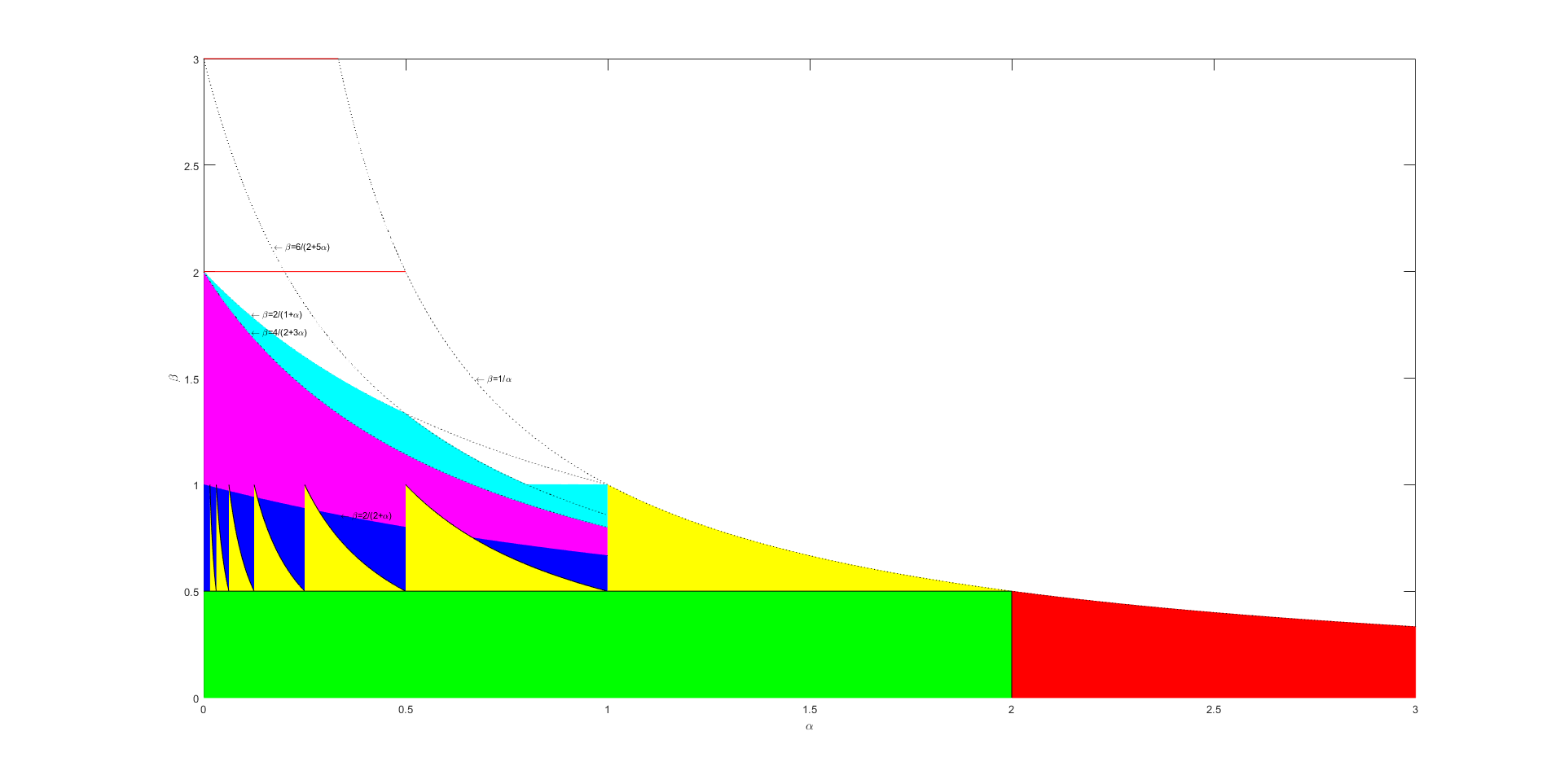}
 \caption{ A sketch of  $\mathcal{F}(g_2)$. The red region contains points $(\alpha, \beta)$ for which $\cG\left(g_2,\alpha, \beta\right)$ is not a frame. All other colors indicate the frame property. The green region is the classical: ``painless expansions" \cite{DaGroMe}, and the yellow region is the result from \cite{Olehon1}. The blue and the magenta regions are respectively from \cite{Lemniel} and \cite{Olehon}. The cyan region is from \cite{AtiKouOko1, AtiKouOko2}.}
\label{fig:figure1}
\end{figure}

\begin{prop}[Frame set of the $2-$spline, $g_2$]\label{frsetg_2} The following statements hold. 
 \begin{enumerate}
\item[(a)]  If $(\alpha, \beta)\in \mathcal{F}(g_2)$, then $\alpha \beta <1$ and $\alpha<2$      \cite{DaGroMe}. This is illustrated by the green region in Figure~\ref{fig:figure1}. 
\item[(b)]  Assume that $1\leq \alpha<2$ and $0<\beta<\frac{1}{\alpha}$. Then, $(\alpha, \beta)\in \mathcal{F}(g_2)$ \cite{Olehon1}. This is illustrated by part of the Yellow region in Figure~\ref{fig:figure1}.
\item[(c)]  Assume that $0<\alpha<2$, and $0<\beta\leq\frac{2}{2+\alpha}$. Then, $(\alpha, \beta)\in \mathcal{F}(g_2)$, and there is a unique 
dual $h\in L^2(\R)\cap L^\infty(\R)$ such that $\text{supp}\, h\subseteq\left[-\frac{\alpha}{2},\frac{\alpha}{2}\right]$ \cite{Olehon}. This is illustrated by the blue region in Figure~\ref{fig:figure1}.
\item[(d)]  Assume that $0<\alpha<2$, and $\frac{2}{2+\alpha}<\beta\leq\frac{4}{2+3\alpha}$. Then, $(\alpha, \beta)\in \mathcal{F}(g_2)$, and there is a unique 
dual $h\in L^2(\R)\cap L^\infty(\R)$ such that $\text{supp}\, h\subseteq\left[-\frac{3\alpha}{2},\frac{3\alpha}{2}\right]$ \cite{Lemniel}. This is illustrated by the magenta region in Figure~\ref{fig:figure1}.
\item[(e)]  Assume that $0<\alpha<1/2$, and $\frac{4}{2+3\alpha}< \beta\leq \frac{2}{1+\alpha}$. Then, $(\alpha, \beta)\in \mathcal{F}(g_2)$, and there is a unique 
dual $h\in L^2(\R)\cap L^\infty(\R)$ such that $\text{supp}\, h\subseteq\left[-\frac{5\alpha}{2},\frac{5\alpha}{2}\right]$  \cite{AtiKouOko2}. This is illustrated by the cyan region in Figure~\ref{fig:figure1}.
\item[(f)]  Assume that $\tfrac{1}{2}\leq \alpha\leq \tfrac{4}{5}$, and $\frac{4}{2+3\alpha}<\beta\leq\frac{6}{2+5\alpha}$, with $\beta>1$. Then, $(\alpha, \beta)\in \mathcal{F}(g_2)$, and there is a unique 
dual $h\in L^2(\R)\cap L^\infty(\R)$ such that $\text{supp}\, h\subseteq\left[-\frac{5\alpha}{2},\frac{5\alpha}{2}\right]$  \cite{AtiKouOko2}. This is illustrated by the cyan region in Figure~\ref{fig:figure1}.
\item[(g)]  Assume that $\tfrac{2}{3}\leq \alpha\leq 1$, and $\frac{4}{2+3\alpha}<\beta < 1$. Then, $(\alpha, \beta)\in \mathcal{F}(g_2)$, and there is a unique compactly supported 
dual $h\in L^2(\R)\cap L^\infty(\R)$ \cite{AtiKouOko1}. This is illustrated by the cyan region in Figure~\ref{fig:figure1}.
\item[(h)]  If $0<\alpha<2$,  $\beta=2, 3, \hdots$, and $\alpha \beta<1$, then, $(\alpha, \beta)\not\in \mathcal{F}(g_2)$  \cite{Grojan}. This is illustrated by the red horizontal lines in Figure~\ref{fig:figure1}.
 \end{enumerate}
\end{prop} 

These results are illustrated  in Figure~\ref{fig:figure1}, where except for the red regions, all other regions are contained in $\mathcal{F}(g_2)$. 
 For the proofs  we refer to \cite{AtiKouOko1, AtiKouOko2, Olehon, Olehon1, chri2003, Grojan, Lemniel, Kloosto}.  But we point out that  the main idea in establishing parts (c--g) is based on the following result due to Janssen \cite{Jan1}. Before stating it we recall that for $\alpha, \beta>0$ and $g\in L^2(\R)$, the Gabor system $\mathcal{G}(g, \alpha, \beta)$ is called  \emph{a Bessel sequence}  if only the upper bound in~\eqref{gabframe} is satisfied for some $B>0$.

\begin{prop}[Sufficient and necessary condition for dual Gabor frames]\cite{Jan1} Let $\alpha, \beta >0$ and $g, h \in L^2(\R).$ The Bessel sequences $\mathcal{G}(g, \alpha, \beta)$ and $\mathcal{G}(h, \alpha, \beta)$ are  dual Gabor frames if and only if 
$$\sum_{k\in \Z}\overline{g(x-n/\beta-k\alpha)}h(x-k\alpha)=\beta \delta_{n,0}\qquad \text{a.e.}\, x \in [0, \alpha].$$
\end{prop}

Using this result with $g=g_N$ and imposing that $h$ is also compactly supported, leads one to seek an appropriate (finite) square matrix from the (infinite) linear system 
$$\sum_{k\in \Z}g_N(x-\tfrac{\ell}{\beta}+k\alpha)h(x+k\alpha)= \beta \delta_{\ell}\qquad \text{for\, almost\, every}\,  x \in [-\tfrac{\alpha}{2}, \tfrac{\alpha}{2}].$$  In particular, \cite{AtiKouOko1} shows that the region $\{(\alpha, \beta) \in \R^2_+: 0<\alpha \beta <1\}$ can be partitioned in subregions $R_m$, $m\geq 1$,  such that a $(2m-1)\times (2m-1)$ matrix $G_m$ can be extracted from the above system leading to $$G_m(x)\begin{bmatrix} h(x+(1-m)\alpha)\\
\vdots\\ h(x)\\ \vdots\\ h(x+(m-1)\alpha)\end{bmatrix}=\begin{bmatrix}0\\ \vdots\\ \beta \\ \vdots\\ 0\end{bmatrix}\qquad \text{for\, almost\, every}\ x\in [-\alpha/2, \alpha/2].$$
Choosing $N=2$  results in parts (c--g) of Proposition~\ref{frsetg_2}, for the cases $m=1,2,$ and $3$. For these cases, one proves that the matrix $G_m(x)$ is invertible for $\text{for\, almost\, every\,} x\in  [-\alpha/2, \alpha/2].$ However, only a subregion for the case $m=3$ has been settled in \cite{AtiKouOko2}. It is also known that the remaining part of this subregion contains some obstruction points, for example the line $\beta=2$ in Figure~\ref{fig:figure1}. Nonetheless, it seems that one should be able to prove that the region $$\{ (\alpha, \beta): \tfrac{1}{2}\leq \alpha <1, \, \, \tfrac{6}{2+5\alpha}\leq \beta < \tfrac{2}{1+\alpha}, \, \, \beta>1\}$$ is also contained in $\mathcal{F}(g_2)$. But this is still open.

We end this section by observing that the frame set problem is a special case of the more general question of characterizing the \emph{full frame set} $\mathcal{F}_{\textrm{full}}(g)$  of a function $g$, where $$\mathcal{F}_{\textrm{full}}(g)=\{ \Lambda \subset \R^2:\, \mathcal{G}(g, \Lambda)\, \textrm{is\, a \, frame}\}$$ where  $\Lambda$ is the lattice $\Lambda=A\Z^2 \subset \R^2$ with $\text{det}A\neq 0$. The only general result known in this case is for $g(x)=e^{-a|x|^2}$ with $a>0$ in which case \cite{Lyuba, Seip} $$\mathcal{F}_{\textrm{full}}(g)=\{\Lambda \subset \R^2:\, \text{Vol}\Lambda <1\},$$where  the volume of $\Lambda$ is defined by $\text{Vol}(\Lambda)= |\text{det} A|$.

 \section{Wilson Bases}\label{sec4}
 By the BLT (Proposition~\ref{blt}) and Proposition~\ref{density_gabf}(iii), we know that $\mathcal{G}(g, \alpha, 1/\alpha)$ cannot be an ONB if $g$ is well-localized in the time-frequency plane.  To overcome the BLT, K.~G.~Wilson introduced an ONB $\{\psi_{n, \ell}, n\in \N_0, \ell \in \Z\}$, where $\psi_{0, \ell}(x)=\psi_\ell(x)$ and for $n\geq 1$, $\psi_{n, \ell}(x)=\psi_\ell(x-n)$, and such that  $\hat{\psi}_{n, \ell}$ is localized around $\pm n$, that is, $\psi_{n, \ell}$ is a bimodal  function. Wilson presented numerical evidences that this system of functions is an ONB for $L^2(\R)$. In 1992, Daubechies, Jaffard, and Journ\'e formalized Wilson's ideas and constructed examples of bimodal Wilson bases generated by smooth functions. To be specific,  the Wilson system associated with a given  function $g\in L^2$,  is $\mathcal{W}(g)=  \{ \psi_{j,m}: j \in \mathbb Z, m \in \mathbb N_0\}$ where

 \begin{align}\label{orig-wils}
 \psi_{j, m}(x)=\begin{cases}g(x-j) \, & \, \text{if} \,  j\in \Z\\
 \tfrac{1}{\sqrt{2}}T_{\tfrac{j}{2}}(M_{m}+(-1)^{j+m}M_{-m})g(x)\, & \text{if} \,   (j,m)\in \Z \times \N,
 \end{cases}
 \end{align}
 which can simply be rewritten as 
 
 \begin{align*}
 \psi_{j, m}(x)=\begin{cases}\sqrt{2}\cos 2\pi m x \ g(x-\tfrac{j}{2}), \, & \, \text{if} \,  j+m \, \text{is\, even}\\
 \sqrt{2}\sin 2\pi m x \ g(x-\tfrac{j}{2}), \, & \, \text{if} \,  j+m \, \text{is\, odd}.
 \end{cases}
 \end{align*}
 
  It is not hard to see $\{\psi_{j,m}\}$ is an ONB for $L^2(\R)$  if and only if  
 \begin{align*}\begin{cases} \|\psi_{j,m}\|=1\,\,  \text{for\, all}\, (j, m)\in \N_0\times \Z\\
\ip{f}{h}= \sum_{j, m}\ip{f}{\psi_{j,m}}\overline{\ip{h}{\psi_{j,m}}}\, \, \text{for\, all}\, f, h\in L^2. 
\end{cases}
\end{align*}
 Assuming that $g$ and $\hat{g}$ are smooth enough, $\hat{g}$ real-valued, one can show that this is equivalent to 
 
 $$\sum_{m\in \Z}\hat{g}(\xi - m)\hat{g}(\xi - m +2j)=\delta_{j,0}\, \text{for\, almost\, every}\, \xi, \, \text{and\, for\, each}\,  j\in \Z.$$ It follows that one can construct compactly supported $\hat{g}$ that will solve this system of  equations. On the other hand, one can convert these equations into a single one by using another time-frequency analysis tool, the Zak transform which we now define. 
 For $f\in L^2(\R)$ we let $Zf: [0, 1)\times [0,1) \to \C$ 
 be given by $$Zf(x, \xi)=\sqrt{2}\sum_{j\in \Z}f(2(x-j))e^{2\pi i j\xi}$$  $Z$ is a unitary map from $L^2(\R)$ onto $L^2( [0,1)^2)$ and enjoys some periodicity-like properties \cite[Chapter 8]{Groc2001}. Using the Zak transform, and  under suitable regularity assumptions on $g$ and $\hat{g}$, one can show that 
 $\{\psi_{j,m}\}$ is an ONB if and only if $$ |Z\hat{g}(x, \xi)|^2 + |Z\hat{g}(x, \xi + \tfrac{1}{2})|^2=2\, \text{for\, almost\, every}\, (x, \xi) \in [0, 1]^2.$$ 
 
 Real-valued functions $g$ solving this equation, can be constructed with the additional requirement that both $g$ and $\hat{g}$ have exponential decay.

 To connect this Wilson system to Gabor frame, we use once again the Zak transform, and observe that the frame operator of the Gabor system $\mathcal{G}(g, 1, 1/2)$ is a multiplication operator in the Zak transform domain, that is $$ZS_g f(x, \xi)=M(x, \xi)Zf(x, \xi)$$ where $M(x, \xi)=|Zg(x, \xi)|^2+|Zg(x, \xi-\tfrac{1}{2})|^2$.  Consequently, $\mathcal{G}(g, 1, 1/2)$ is a tight frame if and only if  $$M(x, \xi)=|Zg(x, \xi)|^2+|Zg(x, \xi-\tfrac{1}{2})|^2=A\, \, \, \text{for\, almost\, every}\, (x, \xi) \in [0, 1]^2,$$ where $A$ is a constant.   These ideas were used in \cite{DaJaJo91} resulting in  the following.

\begin{prop}[\cite{DaJaJo91}]\label{daub-wil} There exist unit-norm real-valued functions $g \in L^{2}(\mathbb R)$  with the property that both $g$ and $\hat{g}$ have exponential decay and such that  the Gabor system $\mathcal{G}(g, 1, 1/2)$ 
is a tight frame for $L^{2}(\R)$ if and only if   the associated  Wilson system
$\mathcal{W}(g)$ 
is an orthonormal basis for $L^2(\R)$.
\end{prop}

Proposition~\ref{daub-wil} also  provides  an alternate view of the Wilson ONB. Indeed,  each function in~\eqref{orig-wils} is a linear combination of at most two Gabor functions  from a tight Gabor frame $\mathcal{G}(g, 1, 1/2)$  of redundancy $2$. Furthermore,  such Gabor systems can be constructed  so as the generators are well-localized in the time-frequency plane.  Suppose now that we are given a tight Gabor system $\mathcal{G}(g, \alpha, \beta)$ where $(\alpha\beta)^{-1}=N\in \N$ where $N> 2$. Hence, the redundancy of this tight frame is $N$. Can a Wilson-type ONB (generated by well-localized window) be constructed from this system by taking appropriate linear combinations? This  problem was posed by Gr\"ochenig  for the case $\alpha=1$ and $\beta=1/3$   \cite[Section 8.5]{Groc2001}, and to the best of our knowledge it is still open.
  If one is willing to give up on the orthogonality, Wojdy{\l}{\l}o \cite{Woj07} proved the existence of  Parseval Wilson-type frames for $L^2(\R)$ from Gabor tight frames of redundancy $3$. More recently, explicit examples have been constructed  starting from Gabor tight frames of redundancy $\tfrac{1}{\beta}\in \N$ where $N\geq 3$.

\begin{prop}\label{Par_Wil}\cite{DBKO18}  For any $\beta \in [1/4,1/2)$ there exists a real-valued function  $g \in S(\R)$ such that the following equivalent statements hold.

\begin{enumerate}
\item[(i)] $\mathcal{G}(g, 1, \beta)$  is a tight Gabor frame of redundancy $\beta^{-1}$. 
\item[(ii)] The associated Wilson system given by 
\begin{equation}\label{wilson-sys}
\mathcal{W}(g, \beta)= \{ \psi_{j,m}: j \in \mathbb Z, m \in \mathbb N_0\}
\end{equation} where 

\begin{align}\label{dks}
\psi_{j,m}(x)=\begin{cases}  \sqrt{2\beta}g_{2j, 0}(x)= \sqrt{2\beta} g(x-2\beta j) &  \text{if} \ \  j\in \mathbb Z, m=0,\\ 
\sqrt{\beta} \left[ e^{-2\pi i \beta j  m}g_{j,m}(x) + (-1)^{j+m}e^{2\pi i \beta j  m} g_{j,-m}(x)\right] & \text{if} \  \  (j,m)\in \mathbb Z \times \mathbb N
\end{cases}
\end{align}
is a Parseval frame for $L^2(\R)$.
\end{enumerate}
If in addition $\beta=\frac{1}{2n}$ where $n$ is any odd natural number, then we can choose $g$ to be real-valued such that both $g$ and $\hat{g}$ have exponential decay. 
\end{prop}

To turn these Parseval (Wilson) systems in  ONBs, one needs to ensure that $\|\psi_{j,m}\|_2=1$ for all $j, m$. This requires in particular that $\|g\|=\tfrac{1}{\sqrt{2\beta}}$, which seems to be incompatible with all the other conditions imposed $g$. It has then been suggested in \cite{DBKO18} that to obtain a Wilson ONB with redundancy different from $2$, one must modify in a fundamental way~\eqref{dks}. For example, if we want to have a Wilson ONB with $\alpha=1, \beta=1/3$, it seems that one should take  linear combinations of three Gabor atoms instead of the two in Proposition~\ref{Par_Wil}. While we have no proof of this claim, it seems to be supported by a recent  construction of multivariate Wilson ONBs which is not  a tensor  products on $1$-Wilson ONBs. This new approach was introduced  in \cite{BJLO17}, where a relationship between these bases and the theory of Generalized Shift Invariant Spaces (GSIS) \cite{RonShen95, RonShen97, RonShen05} was used to construct (non-separable) well-localized Wilson ONB for $L^2(\R^d)$ starting from tight Gabor frames of redundancy $2^k$ where $k=0, 1, 2, \hdots d-1$.  In particular, the functions in the corresponding Wilson systems are linear combinations of $2^k$ elements from the tight Gabor frame.

 \section{HRT}\label{sec5}
 Any application involving Gabor frames, a truncation is needed, and  one considers only a finite number of Gabor atoms. As such, and from a numerical point of view determining the condition number of the projection matrix $$P_{N, K}=\sum_{n=-N}^N\sum_{k=-K}^K \ip{\cdot}{M_{k\beta}T_{n\alpha}g}M_{k\beta}T_{n\alpha}g$$ for $N, K\geq 1$ is useful. In fact, and beyond any numerical considerations, one could ask if this operator is invertible, which will be the case if $\{M_{k\beta}T_{n\alpha}g, \, |n|\leq N, |k|\leq K\}$ was linearly independent. Clearly this is the case if the starting Gabor frame was an ONB. However, and in general this is not known. In fact, this is a special case of a broader problem that we consider in this last section.  This fascinating (due in part to the simplicity of its statement) open problem that was posed in 1990 by C.~Heil, J.~Ramanathan, and P.~Topiwala, and is now referred to as the  HRT conjecture \cite{Heil06, HRT96}.

\begin{conjecture}[The HRT Conjecture]\label{hrt}
Given any $0\neq g \in L^{2}(\R)$  and $\Lambda=\{(a_k, b_k)\}_{k=1}^{N} \subset \R^2$, $\mathcal{G}(g, \Lambda)$ is a linearly independent set in $L^2(\R)$, where $$\mathcal{G}(g, \Lambda)=\{e^{2\pi i b_{k}\cdot}g(\cdot -a_{k}), \, k=1,2, \hdots, N\}.$$ 
\end{conjecture}  
 To be more explicit, the conjecture claims  the following: Given $c_1, c_2, \hdots, c_N\in \C$ such that  
\begin{equation}\label{hrt-explicit}
\sum_{k=1}^Nc_kM_{b_{k}}T_{a_{k}}g(x)=\sum_{k=1}^Nc_ke^{2\pi i b_kx} g(x-a_k)=0\, \, \text{for\, almost\, every}\, x \in \R \, \implies c_1=c_2=\hdots=c_N=0.
\end{equation}
The conjecture is still generally open even  if one assumes that $g\in S(\R)$, the space of $C^\infty$ functions that decay faster than any polynomial.

Observe that  for a given $\Lambda=\{(a_k, b_k)\}_{k=1}^{N} \subset \R^2$, and $g\in L^2(\R)$, we can always assume that  $(a_1, b_1)=(0, 0)$, if not,   applying $M_{-b_{1}}T_{-a_{1}}$ to $\mathcal{G}(g, \Lambda)$ results in $\mathcal{G}(M_{-b_{1}}T_{-a_{1}}g, \Lambda')$ where $\Lambda'$ will include the origin. In addition, by rotating and scaling if necessary, we may also assume that $\Lambda$ contains $(0,1)$. This will result in unitarily changing $g$. Finally, by applying a shear matrix, we may  assume that $\Lambda$ contains $(a, 0)$ for some $a\neq 0$. Consequently, given $\Lambda=\{(a_k, b_k)\}_{k=1}^{N} \subset \R^2$ with $N\geq3$, we shall assume that $\{(0,0), (0,1), (a, 0)\}\subseteq \Lambda$, for some $a\neq 0$. 

To illustrate some of the difficulties arising in investigating this problem, we would like to give some ideas of the proof of the conjecture when $N\leq 3$ and  $0\neq g\in L^2(\R)$. Let us first consider the case $N=2$, and from the above observations we can assume that $\Lambda=\{(0,0), (0,1)\}$. Suppose that $c_1, c_2\in \C$ such that $c_1g+c_2M_1g=0$. This is equivalent to $$(c_1 +c_2e^{2\pi i x})g(x)=0$$ Since $g\neq 0$ and $c_1+c_2e^{2\pi ix}$ is a trigonometric polynomial, we see that $c_1=c_2=0$. 

Now consider the case $N=3$, and assume that $\Lambda=\{(0,0), (0,1), (a, 0)\}$ where $a> 0$  is  such that $\mathcal{G}(g, \Lambda)$ is linearly dependent. Thus there are non-zero complex numbers $c_1,c_2$ such that 
 \begin{equation*}\label{3hrt-base}
 g(x-a)=(c_1 +c_2e^{2\pi i x})g(x)=P(x)g(x)\qquad \forall \, x\in S
 \end{equation*}
where $S\subset \text{supp}(g)\cap(0,1)$ has positive Lebesgue measure. Note that $P(x)$ is a $1$-periodic trigonometric polynomial, that is nonzero almost everywhere. We can now iterate~\eqref{3hrt-base} along $\pm na$ for $n>0$ to obtain

\begin{align*}
\begin{cases}
g(x-na)=g(x)\prod_{j=0}^{n-1}P(x-ja)=g(x)P_{n}(x)\\
g(x+na)=g(x-a)\prod_{j=0}^nP(x+ja)^{-1}=g(x)Q_n(x)
\end{cases}
\end{align*} 
Consequently,  $g(x+na)=g(x)Q_n(x)=g(x)P_{n}(x+na)^{-1}$ implying that 
\begin{equation}\label{q-p}
Q_n(x)=P_n(x+na)^{-1}\qquad  x\in S
\end{equation}
In addition, using the fact that $g\in L^2(\R)$ one can conclude that 
\begin{equation}\label{qpblowup}
\lim_{n\to \infty}P_n(x)=\lim_{n\to \infty}Q_n(x)=0\qquad a.e.\, x\in S
\end{equation}
However, one can show that~\eqref{q-p} and~\eqref{qpblowup} cannot hold simultaneously by distinguishing the case $a\in \Q$ and the case $a$ is irrational. Hence, the HRT conjecture holds when $\#\Lambda=3$. We refer to \cite{HRT96} for details. 

In addition to the fact that the HRT conjecture is true for any set of $3$ distinct points, the known results  generally fall into the following  categories. 
\begin{prop}[HRT for arbitrary set $\Lambda\subset \R^2$]\label{hrt-lambda}
Suppose that $\Lambda\subset \R^2$ is a finite subset of distinct points. Then the HRT conjecture holds in each of the following cases. 
\begin{enumerate}
\item[(a)] $g$ is compactly supported, or just supported within a half-interval $(-\infty, a], $ or $[a, \infty)$ \cite{HRT96}.
\item[(b)] $g(x)=p(x)e^{-\pi x^2}$ where $p$ is a polynomial \cite{HRT96}.
\item[(c)] $g$ is such that $\lim_{x \to \infty}|g(x)|e^{cx^{2}}=0$ for all $c>0$ \cite{BowSpee13}. 
\item[(d)]  $g$ is such that $\lim_{x \to \infty}|g(x)|e^{cx\log x}=0$ for all $c>0$ \cite{BowSpee13}. 
\end{enumerate}
\end{prop}

\begin{figure}[!h]
\begin{center}
 \includegraphics[scale=0.3]{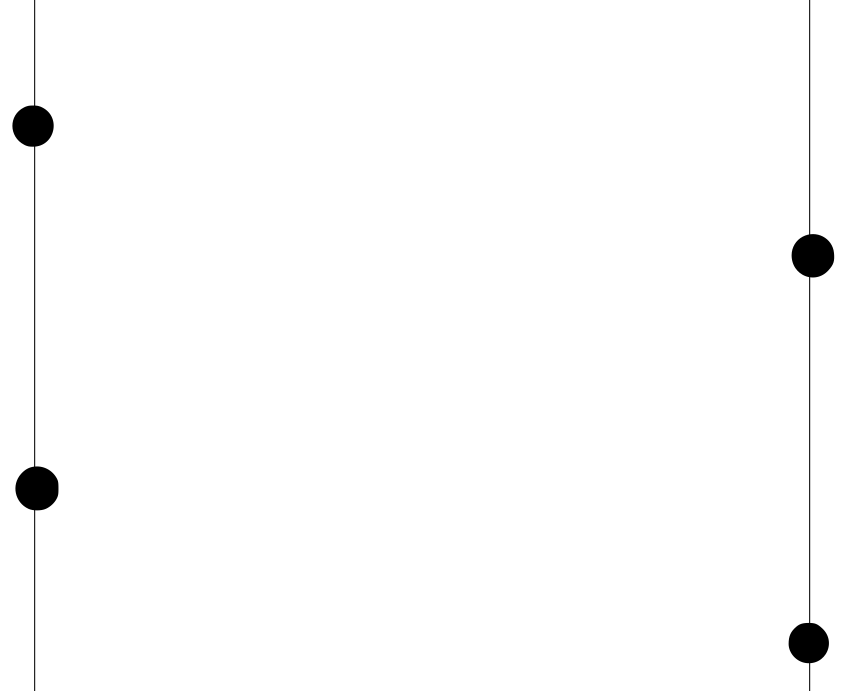}
 \caption{Example of a  $(2,2)$ configuration.}
\label{fig:figure2}
\end{center}
\end{figure}

\begin{prop}[HRT for arbitrary $g\in L^2(\R)$]\label{hrt-l2}
Suppose that  $0\neq g\in L^2(\R)$ is arbitrary.  Then the HRT conjecture holds in each of the following cases. 
\begin{enumerate}
\item[(a)]  $\Lambda$ is a finite set with  $\Lambda \subset A(\Z^2) + z$ where $A$ is a full rank $2\times 2$ matrix and $z\in \R^2$ \cite{Lin99}. In particular, Conjecture~\ref{hrt} holds when $\# \Lambda \leq 3$ \cite{HRT96}.
\item[(b)]  $\# \Lambda =4 $ where two of the four points in $\Lambda$ lie on a line and the remaining two points lie on a second parallel line \cite{Dem10, DemZah12}. Such set $\Lambda$ is called a $(2,2)$ configuration, see Figure~\ref{fig:figure2} for an illustrative example. 
\item[(c)]  $\Lambda$ consists of collinear points \cite{HRT96}.
\item[(d)]  $\Lambda$ consists of $N-1$ collinear and equi-spaced  points, with the last point located off this line \cite{HRT96}.
\end{enumerate}
\end{prop}

We observe that when $\Lambda$ consists of collinear points, the HRT conjecture reduces to the question of linear independence of (finite) translates of $L^2$ functions that was investigated by Rosenblatt \cite{Ros08}. Recently, Currey and Oussa showed that the HRT conjecture is equivalent to the question of linear independence of finite translates of square integrable functions on the Heisenberg group \cite{CurOus18}. 

To date and to the best of our knowledge Proposition~\ref{hrt-lambda} and Proposition~\ref{hrt-l2} are the most general known results on the HRT conjecture. Nonetheless, we give a partial list of known results when one makes restrictions on both the function $g$ and the set $\Lambda$. For an extensive survey on the state of the HRT conjecture we refer to \cite{Heil06, HeiSpe15}, and to \cite{Groc14} for some perspectives from a numerical point of view. 

\begin{figure}[!h]
 \includegraphics[scale=0.3]{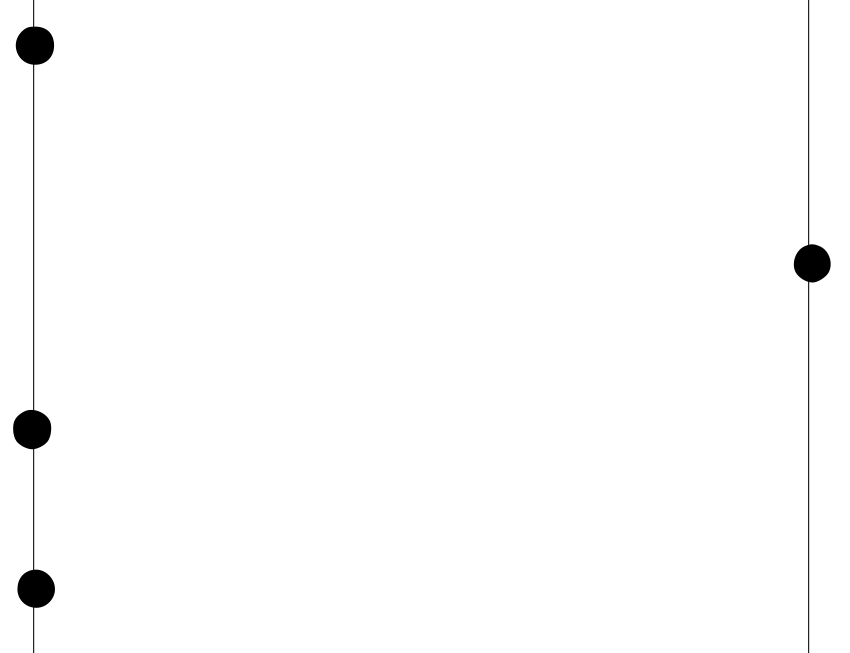}
 \caption{ Example of a $(1,3)$ configuration.}
\label{fig:figure3}
\end{figure}

\begin{prop}[HRT in special cases]\label{hrt-special}
The HRT conjecture holds in each of the following cases. 
\begin{enumerate}
\item[(a)]  $g \in S(\R)$, and  $\# \Lambda =4 $ where three  of the four points in $\Lambda$  lie on a line and the fourth point is off this line  \cite{Dem10}. Such set $\Lambda$ is called a $(1,3)$ configuration,  see Figure~\ref{fig:figure3} for an illustrative example. 
\item[(b)]  $g \in L^2(\R)$ is ultimately positive, and  $\Lambda =\{(a_k, b_k)\}_{k=1}^N \subset \R^2$ is such that $\{b_k\}_{k=1}^N$ are independent over the rationals $\Q$ \cite{BeBo13}. 
\item[(c)]  $\# \Lambda =4 $, when $g \in L^2(\R)$ is ultimately positive,  $g(x)$ and $g(-x)$ are ultimately decreasing \cite{BeBo13}.  
\item[(d)]  $g \in L^2(\R)$ is real-valued, and  $\# \Lambda =4 $  is a $(1,3)$ configuration  \cite{KAO-hrt18}.
\item[(e)]  $g \in S(\R)$ is real-valued functions in $\S(R)$ and $\#\Lambda=4$  \cite{KAO-hrt18}.
\end{enumerate}
\end{prop}

Recently, some of the  techniques used to establish the HRT for $(2,2)$ configurations were extended to deal with some special $(3,2)$ configurations \cite{KAO-hrt18}. From these results, and when restricting to real-valued functions, it was concluded that the HRT holds for certain sets of four points. We briefly describe this method here.

Let  $\Lambda=\{(0,0), (0,1), (a_0, 0), (a, b)\}$, and assume that $\Lambda$ is neither a  $(1,3)$, nor a $(2,2)$ configuration. Let $0\neq g\in L^2(\R)$ be a real-valued function. Suppose that $\mathcal{G}(g, \Lambda)$ is linearly dependent. Then there exist $0\neq c_k \in \C$, $k=1,2,3$ such that $$T_{a_{0}}g=c_1g + c_2M_1g+c_3M_bT_ag$$ Taking the complex conjugate of this equation leads to $$T_{a_{0}}g=\overline{c}_1g+\overline{c}_2M_{-1}g+\overline{c}_3M_{-b}T_ag$$ Taking the difference of these two equations gives $$ (c_1-\overline{c}_1)g+ c_2M_1g-\overline{c}_2g + c_3M_bT_ag-\overline{c}_3M_{-b}T_ag=0.$$ Since $c_2, c_3\neq 0$ we conclude that $\mathcal{G}(g, \Lambda')$ where $\Lambda'=\{(0,0),(0,1), (0, -1), (a,b), (a, -b)\}$ is a (symmetric) $(3,2)$ configuration,  is linearly dependent.  Consequently, we have proved the following result.

\begin{prop}\label{conf2-3} Let $0\neq g \in L^2(\R)$ be a real-valued function. Suppose that $(a, b)\in \R^2$ is such that $\mathcal{G}(g, \Lambda_0)$ is linearly independent where $\Lambda_0=\{(0,0), (0,1), (0, -1), (a, b), (a, -b)\}$. Then for all $0\neq c\in \R$, $\mathcal{G}(g, \Lambda)$ is linearly independent where $\Lambda=\{(0,0), (0,1), (c, 0), (a, b)\}$.
\end{prop} 

In \cite[Theorem 6, Theorem 7]{KAO-hrt18} it was proved that the hypothesis of Proposition~\ref{conf2-3} is satisfied  when $g\in L^2(\R)$ (not necessarily real-valued) for certain values of $a,$ and $b$. These results were viewed as a  \emph{restriction principle} for the HRT, whereby proving the conjecture for special sets of $N+1$ points one can establish it for certain related sets of $N$ points. In addition, a related \emph{extension principle}, which can be viewed as an induction-like technique was introduced.  The premise of this  principle is based on the following question.  Suppose that  the HRT conjecture holds for all $g \in L^{2}(\R)$ and a set   $\Lambda=\{(a_k, b_k)\}_{k=1}^{N} \subset \R^2$. For which points $(a, b) \in \R^2 \setminus \Lambda$ will the conjecture remain true for the same function $g$ and the new set $\Lambda'=\Lambda \cup \{(a,b)\}$?  

We  elaborate on  this method  for $\#\Lambda=3$. Let $g\in L^2(\R)$ with $\|g\|_2=1$ and suppose that $\Lambda=\{(0,0), (0,1), (a_0, 0)\}$. We denote $\Lambda'=\Lambda\cup\{(a, b)\}= \{(0,0), (0,1), (a_0, 0), (a, b)\}$. Since, $\mathcal{G}(g, \Lambda)$ is linearly independent, the Gramian of this set of function is a positive definite matrix. We recall that the Gramian of a set of $N$ vectors $\{f_k\}_{k=1}^N\subset L^2(\R)$ is the (positive semi-definite $N\times N$ matrix $(\ip{f_k}{f_\ell})_{k, \ell =1}^N$. In the  case at hand,  the $4\times 4$ Gramian matrix $G:=G_{g}(a,b)$  of $\mathcal{G}(g, \Lambda')$ can be written in the following block structure:

\begin{equation}\label{gram3}
G=\begin{bmatrix}A & u(a,b)\\ u(a,b)^* & 1\end{bmatrix}
\end{equation}

 where $A$ is the $3\times 3$ Gramian of $\mathcal{G}(g, \Lambda)$ and 
 
 \begin{equation*}
u(a,b)=\begin{bmatrix}V_{g}g(a, b)\\ V_gg(a, b-1)\\ e^{-2\pi i a_{0}b}V_gg(a-a_0, b)\end{bmatrix}
\end{equation*} and 
$ u(a,b)^{*}$ is the adjoint of  $u(a,b)$.  
By construction $G$ is positive semi-definite for all $(a, b)\in \R^2$ and we seek the set of points $(a, b)\in \R^2\setminus \Lambda$ such that $G$ is positive definite. We can encode this information into the determinant of this matrix, or into a related function $F: \R^2 \to [0, \infty)$   given by 
\begin{equation}\label{maintool}
F(a,b)=\ip{A^{-1}u(a,b)}{u(a,b)}.
\end{equation}
 
 The following  was proved in \cite{KAO-hrt18}.
 
 \begin{prop}[The HRT Extension function]\label{ext-kao}
 Given the above notations the function $F$ satisfies the following properties. 

 \begin{enumerate}
\item[(i)] $0\leq F(a,b)\leq 1$ for all $(a, b) \in \R^2$, and moreover, $F(a, b)=1$ if $(a, b)\in \Lambda.$ 
\item[(ii)] $F$ is uniformly continuous and $\lim_{|(a,b)|\to \infty}F(a,b)=0$.
\item[(iii)] $\iint_{\R^{2}}F(a,b)da db=3$.
\item[(iv)] $\det{G_{g}(a,b)}=(1-F(a,b))\det{A}$.
\end{enumerate}
 
 Consequently, there exists $R>0$ such that the HRT conjecture holds for $g$ and $\Lambda'= \Lambda\cup\{(a, b)\}= \{(0,0), (0,1), (a_0, 0), (a, b)\}$ whenever $|(a, b)|>R.$
 \end{prop}
 
%This results suggest that extending the HRT conjecture one point at the time is a  ``local'' problem. It can also be used to numerically  investigate the conjecture. 
%

We conclude the paper by elaborating on the case $\Lambda=4$.  Let $\Lambda\subset \R^2$ contain $4$ distinct points, and  assume without loss of generality that  $\Lambda=\{(0,0), (0,1), (a_0, 0), (a, b)\}$. 
 
 When $b=0$ and $a=-a_0$ or $a=2a_0$,  then  $\Lambda$ is a $(1,3)$ configuration with the additional fact that its three collinear points are equi-spaced. This case is  handled by Fourier methods as was done in \cite{HRT96}; see Proposition~\ref{hrt-l2} (d).  But, for general $(1,3)$ configurations, the Fourier methods are  ineffective. Nonetheless, this case was considered by Demeter \cite{Dem10}, who proved that the HRT conjecture holds for all $(1,3)$ configurations when  $g\in \S(\R)$, and for a family of $(1,3)$ configurations when $g\in L^2(\R)$.  It was latter proved by Liu that in fact, the HRT holds for all functions  $g\in L^2(\R)$ and for almost all (in the sense of Lebesgue measure) $(1,3)$ configurations  \cite{wliu16}. In fact, more is true, in the sense that for $g\in L^2(\R)$, there exists at most one (equivalence class of) $(1, 3)$ configuration $\Lambda_0 $ such that $\mathcal{G}(g, \Lambda_0)$ is linearly dependent \cite{KAO-hrt18}. Here, we say that two sets $\Lambda_1$ and $\Lambda_2$ are equivalent if there exists a symplectic matrix $A\in SL(2, \R)$ (the determinant of $A$ is $1$) such that $\Lambda_2=A\Lambda_1$.  However, it still not known if the HRT holds for all $(1,3)$ configurations when $g\in L^2$. 
 
Next if $b=1$ with $a\not\in \{ 0, a_0\}$, or if $a=a_0$ with $b\neq 0$ then $\Lambda$ is a $(2,2)$ configuration. Demeter and Zaharescu \cite{DemZah12} established the HRT for all $g\in L^2$ and all such configurations. % In dealing with both $(1,3)$ and $(2,2)$ configurations, some delicate estimates and number theoretical properties of $a$ and $b$ were exploited. 
 
 Consequently, to establish the HRT conjecture for all sets of four distinct points and all $L^2$ function, one needs to focus on  \newline
\noindent $\bullet$ showing that there is no equivalence class of  $(1,3)$ configurations for which the HRT fails; and \newline
\noindent $\bullet$ proving the HRT for sets of four points that are neither $(1,3)$ configurations nor $(2,2)$ configurations.  

 For illustrative purposes we pose the following question.
 
 \begin{question}\label{special4} Let $0\neq g\in L^2(\R)$. Prove that  $\mathcal{G} (G, \Lambda)$  is linearly independent in each of the following cases
 \begin{enumerate}
 \item[(a)] $\Lambda=\{(0,0), (0,1), (1, 0), (\sqrt{2}, \sqrt{2})\}$, see \cite[Conjecture 9.2]{Heil06}.
 \item[(b)] $\Lambda=\{(0,0), (0,1), (1, 0), (\sqrt{2}, \sqrt{3})\}$.
 \end{enumerate}  
 
 To be more explicit, the question is to prove that  each of the following two sets are linearly independent

$$\{g(x), g(x-1), e^{2\pi i x}g(x), e^{2\pi i \sqrt{2} x} g(x-\sqrt{2})\}$$ and 

 $$\{g(x), g(x-1), e^{2\pi i x}g(x), e^{2\pi i \sqrt{3} x} g(x-\sqrt{2})\}$$
 
 \end{question}
When $g$ is real-valued, then  part (a) was  proved in \cite{KAO-hrt18}, but nothing can be said for part (b). On the other hand, \cite[Theorem 7]{KAO-hrt18} establishes part (b) when $g\in \S(\R)$.

\section*{Acknowledgment} 
The author acknowledges G.~Atindehou and F.~Ndjakou Njeunje's help in generating the graph included in the paper.

\bibliographystyle{amsplain}
\bibliography{InvGaborA_Bib}

 \end{document}